\documentclass{amsart}
\begin{document}
\title[transcendental numbers]{Transcendental numbers as solutions to arithmetic differential equations}
\author{Alexandru Buium}

\def \Rp{R_p}
\def \Rpi{R_{\pi}}
\def \dpi{\d_{\pi}}
\def \bT{{\bf T}}
\def \cI{{\mathcal I}}
\def \cH{{\mathcal H}}
\def \cJ{{\mathcal J}}
\def \ZN{\bZ[1/N,\zeta_N]}
\def \tA{\tilde{A}}
\def \o{\omega}
\def \tB{\tilde{B}}
\def \tC{\tilde{C}}
\def \alph{A}
\def \bet{B}
\def \bsigma{\bar{\sigma}}
\def \y{^{\infty}}
\def \Ra{\Rightarrow}
\def \uBS{\overline{BS}}
\def \lBS{\underline{BS}}
\def \lB{\underline{B}}
\def \<{\langle}
\def \>{\rangle}
\def \hL{\hat{L}}
\def \cU{\mathcal U}
\def \cF{\mathcal F}
\def \S{\Sigma}
\def \st{\stackrel}
\def \sd{Spec_{\d}\ }
\def \pd{Proj_{\d}\ }
\def \s{\sigma_2}
\def \i{\sigma_1}
\def \bs{\bigskip}
\def \cD{\mathcal D}
\def \cC{\mathcal C}
\def \cT{\mathcal T}
\def \cK{\mathcal K}
\def \cX{\mathcal X}
\def \sX{X_{set}}
\def \cY{\mathcal Y}
\def \cS{X}
\def \cR{\mathcal R}
\def \cE{\mathcal E}
\def \tcE{\tilde{\mathcal E}}
\def \cP{\mathcal P}
\def \cA{\mathcal A}
\def \cV{\mathcal V}
\def \cM{\mathcal M}
\def \cL{\mathcal L}
\def \cN{\mathcal N}
\def \tcM{\tilde{\mathcal M}}
\def \caS{\mathcal S}
\def \cG{\mathcal G}
\def \cB{\mathcal B}
\def \tG{\tilde{G}}
\def \cF{\mathcal F}
\def \h{\hat{\ }}
\def \hp{\hat{\ }}
\def \tS{\tilde{S}}
\def \tP{\tilde{P}}
\def \tA{\tilde{A}}
\def \tX{\tilde{X}}
\def \tcS{\tilde{X}}
\def \tT{\tilde{T}}
\def \tE{\tilde{E}}
\def \tV{\tilde{V}}
\def \tC{\tilde{C}}
\def \tI{\tilde{I}}
\def \tU{\tilde{U}}
\def \tG{\tilde{G}}
\def \tu{\tilde{u}}
\def \chu{\check{u}}
\def \tx{\tilde{x}}
\def \tL{\tilde{L}}
\def \tY{\tilde{Y}}
\def \d{\delta}
\def \e{\chi}
\def \bW{\mathbb W}
\def \bV{{\mathbb V}}
\def \bF{{\bf F}}
\def \bE{{\bf E}}
\def \bC{{\bf C}}
\def \bO{{\bf O}}
\def \bR{{\bf R}}
\def \bA{{\bf A}}
\def \bB{{\bf B}}
\def \cO{\mathcal O}
\def \ra{\rightarrow}
\def \bx{{\bf x}}
\def \f{{\bf f}}
\def \bX{{\bf X}}
\def \bH{{\bf H}}
\def \bS{{\bf S}}
\def \bF{{\bf F}}
\def \bN{{\bf N}}
\def \bK{{\bf K}}
\def \bE{{\bf E}}
\def \bB{{\bf B}}
\def \bQ{{\bf Q}}
\def \bd{{\bf d}}
\def \bY{{\bf Y}}
\def \bU{{\bf U}}
\def \bL{{\bf L}}
\def \bQ{{\bf Q}}
\def \bP{{\bf P}}
\def \bR{{\bf R}}
\def \bC{{\bf C}}
\def \bD{{\bf D}}
\def \bM{{\bf M}}
\def \bZ{{\mathbb Z}}
\def \xtoleqr{x^{(\leq r)}}
\def \hU{\hat{U}}
\def \k{\kappa}
\def \ee{\overline{p^{\k}}}

\newtheorem{THM}{{\!}}[section]
\newtheorem{THMX}{{\!}}
\renewcommand{\theTHMX}{}
\newtheorem{theorem}{Theorem}[section]
\newtheorem{corollary}[theorem]{Corollary}
\newtheorem{lemma}[theorem]{Lemma}
\newtheorem{proposition}[theorem]{Proposition}
\theoremstyle{definition}
\newtheorem{definition}[theorem]{Definition}
\theoremstyle{remark}
\newtheorem{remark}[theorem]{Remark}
\newtheorem{example}[theorem]{\bf Example}
\numberwithin{equation}{section}
%\subjclass[2000]{....}
\address{University of New Mexico \\ Albuquerque, NM 87131}
\email{buium@math.unm.edu} 
%\subjclass[2000]{11E57,12H05}

\maketitle

%\centerline{\it Dedicated to Vasile Brinzanescu on the occasion of his $70$th birthday}

\begin{abstract}
Arithmetic differential equations are analogues of  algebraic differential equations  in which derivative operators acting on functions are replaced by Fermat quotient operators acting on numbers. Now, various remarkable transcendental functions are solutions to algebraic differential equations; in this note  we show that, in a similar way,  some remarkable transcendental numbers (including certain ``periods") are solutions to arithmetic differential equations. Inspired by a recent paper of Manin, we  then speculate on   the possibility of understanding the algebraic relations among periods via Galois groups of arithmetic differential equations. \end{abstract}

\section{Arithmetic differential equations \cite{char, book}}

In the classical theory of differential equations (in the complex $z$-plane ${\mathbb C}$) one starts 
with a ring ${\mathcal R}$ of  holomorphic functions in $z$ that is stable under the derivation operator $\d_z=\frac{d}{dz}$. Then, for $N\geq 1$,  one defines a {\it differential polynomial function} (or simply a {\it $\d$-function}) of order $r$ (in $N$ variables) as 
a set theoretic map $f:{\mathcal R}^N\ra {\mathcal R}$ with the property that there exists a polynomial in $(r+1)N$ variables such that for all $u\in {\mathcal R}^N$ we have 
$f(u)=F(u,\d_zu,...,\d^r_zu)$.
More generally, if $X\subset {\mathbb A}^N$ is a smooth closed subscheme of the affine space, a function $f:X({\mathcal R})\ra {\mathcal R}$ on the set of ${\mathcal R}$-points of $X$ can be called a {\it $\d$-function} of order $r$ if it is the restriction of a $\d$-function of order $r$ on ${\mathcal R}^N={\mathbb A}^N({\mathcal R})$. One can globalize this notion to the case when $X$ is not necessarily affine.
 An {\it algebraic differential equation} is then an equation of the form $f(P)=0$ with $f$ as above and solutions $P\in X({\mathcal R})$. When investigating transcendence properties  one  also fixes a subring $\cO\subset {\mathcal R}$
that is stable under the operator $\d_z$ and one asks that the $\d$-functions under consideration be defined over $\cO$, i.e. the polynomials $F$ have coefficients in $\cO$; a typical example of choice for $\cO$ is, of course,  $\cO={\mathbb C}[z]$, the ring of polynomials. Then one can investigate the algebraic dependence relations over $\cO$ among solutions of various algebraic differential equations. In case the equations are linear, or more generally attached to algebraic groups, this problem is directly related to the differential Galois theory of these equations \cite{kolchin}.

Let us now recall from \cite{char,book} the basics of a theory that   replaces functions by numbers and derivations $\d_z$ by Fermat quotient operators $\d_p$. 
The analogue of the ring ${\mathcal R}$ will be the ring $R$
described as follows. We start with the ring of integers $\bZ$, we let $\bZ_p$ be the $p$-adic completion of $\bZ$, we let $\bZ_p^{ur}$ be the maximum unramified extension of $\bZ_p$ (obtained by adjoining to $\bZ_p$ all the roots of unity of order prime to $p$), and we let $R$ be the $p$-adic completion of $\bZ_p^{ur}$. The ring $R$ plays a central role in number theory and turns out to have other descriptions (e.g. it is the $p$-typical Witt ring $W(k)$ on the algebraic closure $k=\overline{\mathbb F}_p$ of the field ${\mathbb F}_p$ with $p$ elements; it is also the unique complete  discrete valuation ring with maximal ideal generated by $p$ and residue field $k$). It is well known that any element of $R$ can be written uniquely as a series $\sum c_ip^i$ where $c_i$ are either roots of unity of order prime to $p$ or $0$. Also there is a unique ring automorphism $\phi_p:R\ra R$ whose reduction mod $p$ is the Frobenius automorphism of $k$; $\phi_p$ is given by $\phi_p(\sum c_ip^i)=\sum c_i^pp^i$. We then define the {\it Fermat quotient operator} $\d=\d_p:R\ra R$ by the formula $\d_p u=\frac{\phi_p(u)-u^p}{p}$ and we view this operator as an analogue of a derivation. (Note however that $\d$ is {\it not} additive!) Next recall that by a restricted power series with coefficients in $R$ one understands a power series whose coefficients tend to $0$ $p$-adically. Here is one of our main definitions in \cite{char}: a
 set theoretic function $f:R^N\ra R$ is called a {\it $\d$-function} of order $r$ if there exists  a restricted power series $F$ with coefficients in $R$, in $(r+1)N$ variables, such that for all $u\in R^N$, we have
 $f(u)=F(u,\d_p u,...,\d_p^r u)$.
 More generally if $X\subset {\mathbb A}^N$ is a smooth closed subscheme of the affine space then a function $f:X(R)\ra R$ on the set of $R$-points of $X$ is called a {\it $\d$-function} of order $r$ if it is the restriction of a $\d$-function of order $r$ on $R^N={\mathbb A}^N(R)$. 
 If $X$ is an arbitrary (not necessarily affine) smooth scheme over $R$ a function $f:X(R)\ra R$ is called a {\it $\d$-function} of order $r$ if its restriction to any affine open set of $X$ has this property.
 Finally 
 an {\it arithmetic differential equation} on an $X$ as above is, by definition, an equation of the form 
 $f(P)=0$
 where  the solutions $P$ are in $X(R)$. 
 Again, if one wants to investigate transcendence properties of solutions, one fixes a subring $\cO\subset R$ stable under $\d_p$ and one assumes that the $\d$-functions under consideration are ``defined" over $\cO$. The latter may mean a number of things. One can ask, for instance, that the coefficients of the series $F$ defining $f$ belong to $\cO$; this is a natural condition only if $\cO$ is $p$-adically complete and the natural choice, in this case, is $\cO=\bZ_p$. Other meanings of ``defined over $\cO$"  may make sense in special situations, however; for instance one can consider the case (which effectively occurs in interesting examples)  when the series $F$  are algebraic over a polynomial ring with coefficients in $\cO$; in this case it is natural to take $\cO$ not necessarily $p$-adically complete.

 \section{Transcendental functions/numbers appearing as solutions}
 
 \subsection{Exponential function}
 The simplest example we would like consider is that of the algebraic differential equation
 \begin{equation}
 \label{exponential}
 \frac{\d_z u}{u}=\beta
 \end{equation}
 where $\beta\in {\mathcal R}$ is a fixed function and the solutions are $u\in {\mathcal R}^{\times}$ . We may interpret this as an algebraic differential equation 
 on the hyperbola $H\subset {\mathbb A}^2$  defined by $xy-1=0$,
 by viewing  $P:=(u,\frac{1}{u})\in H({\mathcal R})$ as a point of this hyperbola.
 Of course this hyperbola identifies, via the first projection  with  the multiplicative group scheme
 ${\mathbb G}_m$.
 The solutions to this equation are of the form
 \begin{equation}
 \label{expo}
 u=c\cdot \exp\left(\int \beta dz\right),
 \end{equation}
 where $c\in {\mathbb C}$.
 By the way, the problem of understanding the algebraic relations among  these solutions is closely related to the {\it  Ax-Lindemann-Weierstrass} theorem \cite{ax}
 stating that if $\alpha_1,...,\alpha_n$
 are holomorphic functions that
 are algebraic over ${\mathbb C}(z)$ and  are ${\mathbb Q}$-linearly independent modulo ${\mathbb C}$ then $\exp(\alpha_1),...,\exp(\alpha_n)$ are algebraically independent over ${\mathbb C}(z)$. We recall that the latter is the functional analogue of the {\it Lindemann-Weierstrass} theorem asserting that if  $\alpha_1,...,\alpha_n\in {\mathbb C}$ are algebraic over ${\mathbb Q}$ and are linearly independent over ${\mathbb Q}$ then $\exp(\alpha_1),...,\exp(\alpha_n)$ are algebraically independent over ${\mathbb Q}$.
 
 Here is an arithmetic analogue \cite{char, adel3} of the above situation.
 Consider the hyperbola $H={\mathbb G}_m$ over $R$ and the group homomorphism $\psi:H(R)=R^{\times}\ra R$ defined by
\begin{equation}
\label{psimap}
u\mapsto \psi(u)=\frac{1}{p}\log_p\left(\frac{\phi_p(u)}{u^p}\right)=\sum_{n=1}^{\infty}(-1)^{n-1}
\frac{p^{n-1}}{n}\left(\frac{\d u}{u^p}\right)^n
\end{equation}
where $\log_p$ is the $p$-adic logarithm. Then, for any fixed $\beta\in R$, one can consider the arithmetic differential equation
\begin{equation}
\label{e}
\psi(u)=\beta.
\end{equation}
 Note that the ``defining" restricted power series of $\psi$ is {\it not} algebraic over the ring of polynomials with coefficients in $R$. Since $\psi$ in \ref{psimap} is a homomorphism it can be viewed as an analogue of the logarithmic derivative map $u\mapsto \frac{\d_z u}{u}$; hence \ref{e} can be viewed as  an analogue of \ref{exponential}. 
By the way note that \ref{e} is equivalent to each of the following  equations
\begin{equation}\label{dodododo}
\phi_p (u)=\epsilon\cdot u^p\ \ \ \ \text{or}\ \ \ \ 
\d_p u=\alpha\cdot u^p
\end{equation}
where $\epsilon=\exp_p(p\beta)=1+p\alpha$, $\exp_p$ the $p$-adic exponential; in  equation \ref{dodododo} the ``defining" restricted power series  is actually a polynomial. In any case the set of solutions to \ref{dodododo} (and hence to \ref{dodododo} or \ref{e}) consists of all $u\in R^{\times}$ of the form
\begin{equation}
\label{research}
u=\zeta \cdot \exp_p\left( \sum_{n=1}^{\infty} p^n \phi_p^{-n}(\beta)\right)
\end{equation}
where $\zeta\in R^{\times}$, $\d_p \zeta=0$. 
So \ref{research} can be viewed as an analogue of \ref{expo}.
It is worth mentioning that equation \ref{dodododo} is {\it not} an instance of a difference equation in the sense of \cite{SVdP}; indeed difference equations for $\phi_p$ take the form 
\begin{equation}
\label{difference}
\phi_p(u)=\epsilon\cdot u\end{equation}
 rather than \ref{dodododo}. Nevertheless difference equations \ref{difference} are also  examples of arithmetic differential equations and their Galois theory as arithmetic differential equations \cite{book} is quite different from their Galois theory as difference equations \cite{SVdP}.
In terms of transcendence properties one can easily check the following:
there exists  a subset $\Omega$  of the first category in the metric space
$X:=1+pR$
such that for any $u \in X\backslash \Omega$ there exists a subring $\cO$ of $R$, stable under $\d_p$ and containing $\alpha:=\frac{\d_p u}{u^p}$, with the property that  $u$ is transcendental over $\cO$; cf. \cite{adel3}. On the other hand if $u$ satisfies
\ref{dodododo} and lies in $\bZ_p^{ur}$ then it is an easy exercise to see that $u$ is algebraic over any subring $\cO\subset R$ that is stable under $\d_p$ and contains $\alpha$; cf. \cite{adel3}. On a related note, one has, of course,  the profound results due to Mahler, Gelfond, Nesterenko, Philippon, and many others on transcendence (over ${\mathbb Q}$!)  of exponentials in the $p$-adic setting; one prototype of these is Mahler's theorem stating that  if $a\in {\mathbb C}_p$ is algebraic over ${\mathbb Q}$ and $0<|\alpha|_p<p^{-1/(p-1)}$ then $\exp_p(\alpha)$ is transcendental over ${\mathbb Q}$.

The above example has a matrix analogue \cite{adel3} which we briefly mention here.
In the classical theory of differential equations the matrix analogue of \ref{exponential} is a  linear differential equation
\begin{equation}
\label{linear}
\d_z u=\beta \cdot u,
\end{equation}
where $\beta\in {\mathfrak g}{\mathfrak l}_n({\mathcal R})$ is fixed and the solutions are  $u\in GL_n({\mathcal R})$; here ${\mathfrak g}{\mathfrak l}_n$ means $n\times n$ matrices and $GL_n$ means invertible matrices. The simplest arithmetic analogue of \ref{linear} considered in \cite{adel3} is
\begin{equation}
\label{lineara}
\d_p u=\beta \cdot u^{(p)},
\end{equation}
where $\beta\in {\mathfrak g}{\mathfrak l}_n(R)$ is fixed  and the solutions are $u\in GL_n(R)$; here for $u=(u_{ij})$ we wrote $u^{(p)}:=(u_{ij}^p)$. 
Since the $\d$-functions appearing in \ref{lineara} are actually polynomials  
 it is natural to consider algebraic relations among the entries of solutions to \ref{lineara} over non-complete rings $\cO\subset R$ such that $\beta\in {\mathfrak g}{\mathfrak l}_n(\cO)$; this was explored in \cite{adel3}. 
The equation \ref{linear} corresponds to the trivial ``arithmetic connection" on $GL_n$; more generally one can consider ``arithmetic connections"  attached to symmetric/antisymmetric matrices $q\in GL_n(R)$; they lead to equations
\cite{adel2, adel3} of the form
\begin{equation}
\label{linearaa}
\d_p u=\beta \cdot \Phi_p(u),
\end{equation}
where $\Phi_p$ are appropriate matrices of restricted power series depending on $q$ and compatible with $q$ in a way that is reminiscent of the way Chern connections in complex geometry are compatible with hermitian metrics; we refer to \cite{adel2, adel3} for details. Although $\Phi_p$ are power series (rather than polynomials, as in \ref{lineara})  they are ``algebraic" in a certain precise sense
and hence it is still natural to consider algebraic relations among the entries of solutions to \ref{linearaa} over non-complete rings $\cO\subset R$; this is not covered by \cite{adel3} and deserves being investigated.

\subsection{Elliptic functions}
The next example is taken from \cite{manpain,char, BYM}. A remarkable class of differential equations in the complex plane is that of Painlev\'{e} equations, in particular of the Painlev\'{e} VI family which we now describe following the classical work of Fuchs and the modern interpretation of Manin \cite{manpain}. 
Let $E$ be the (smooth projective) elliptic curve over ${\mathbb Q}(z)$ defined by the equation
$y^2=x(x-1)(x-z)$. Then by Fuchs and Manin \cite{manin63},  there is a (unique up to a multiplicative constant) compatible system of $\d$-functions of order $2$, $\psi:E({\mathcal R})\ra {\mathcal R}$, that  are group homomorphisms, where  ${\mathbb C}(z)\subset {\mathcal R}$.  If $P=(X,Y)\in E({\mathcal R})$ is a point of $E$ with $X,Y\in {\mathcal R}$ 
(so $x(P)=X$, $y(P)=Y$)
then
$$
\psi(P) = z(1-z)\left(z(1-z)\d_z^2+(1-2z)\d_z-\frac{1}{4}\right)
\int_{\infty}^{(X,Y)}\frac{dx}{y}=F(X,Y,\d_zX,\d_z^2 X),$$
where $F$ is a rational function with coefficients in ${\mathbb Q}(z)$ in $4$ variables.
(Manin's original notation for $\psi$ above is $\mu$; cf. \cite{manin63}.)
The expression involving the integral is well defined because the linear differential operator in front of the integral annihilates any function of $z$ of the form
$\int_{\gamma}\frac{dx}{y}$ where $\gamma$ is an integral $1$- cycle.
 Functions of the form $\int_{\gamma}\frac{dx}{y}$ are called in \cite{manin} {\it periods-functions}; the various values of periods-functions at special $z$s are then referred to  as  {\it numerical periods}, or simply {\it periods}.
According to \cite{manpain} the Painlev\'{e} VI equation(s) are then of the form
\begin{equation}
\label{pain}
\psi(P)=\sum_{i=0}^3 \lambda_i \cdot y(P+P_i)
\end{equation}
where $P\in E({\mathcal R})$, $\lambda_i\in {\mathbb C}$, $P_0=0$,   $\{P_1,P_2,P_3\}$ are the points of $E({\mathcal R})$ of order $2$, and $P+P_i$ is the sum in the group law of the elliptic curve. The solutions to \ref{pain}, when transcendental,  are the famous {\it Painlev\'{e} (VI) transcendents}  and their consideration was motivated, in the work of Painlev\'{e}, by his search of ``new" transcendental functions, solving differential equations, and having  ``no movable singularity" (by which it is roughly meant that the non-algebraic singularities of the solutions do not depend on the integration constants). There is a vast literature on Painlev\'{e} equations which seem to pop up is a number of basic, but apparently unrelated, mathematical (and physical) contexts; in particular these equations have a remarkable ``Hamiltonian structure". Now, an arithmetic analogue of Manin's $\psi$ was introduced in \cite{char} where it was proved that, for any elliptic curve $E$ over $R$, 
there exists a non-zero $\d$-function $\psi_p:E(R)\ra R$ of order $2$ which is a group homomorphism; this $\psi_p$ is actually unique up to a multiplicative constant, provided $E$  is not a {\it canonical lift} (in particular if $E$ is  without complex multiplication).
Then the equation \ref{pain} makes sense in the arithmetic case by replacing $\psi$ with $\psi_p$, taking  $P\in E(R)$, and $\lambda_i\in R$; more generally one can consider, as analogue of \ref{pain}, the equation
\begin{equation}
\label{suffering}
\psi_p(P)=\sum_{i=0}^3 \lambda_i \cdot y(P+P_i)^{\phi_p^{\nu}}
\end{equation}
where $\nu\in \bZ_{\geq 0}$.
For $\nu=1$ the equation \ref{suffering} was shown in \cite{BYM} to possess a structure reminiscent of Hamiltonian structure; on the other hand
 the solutions  of \ref{suffering} can be viewed as arithmetic analogues of the Painlev\'{e} transcendents, and could be referred to as {\it Painlev\'{e} numbers}. A question immediately invites itself, although, as stated, it is definitely too vague: do Painlev\'{e} numbers have a property that can be viewed as an arithmetic analogue of ``no movable singularity"? Another question can be asked, however, that is quite precise: what are the algebraic relations among solutions of \ref{pain} in the arithmetic case? (In the classical case of differential equations there is ample work on the corresponding question, both classical and modern, cf. \cite{NP} and the bibliography therein.) The degenerate case $\lambda_0=...=\lambda_3=0$ is worth examining separately. In this case it is known (a special case of Manin's theorem of the kernel \cite{manin63}) that if ${\mathcal R}$ is a field algebraic over ${\mathbb C}(t)$
 then any solution $P$ to \ref{pain} (i.e. any $P\in E({\mathcal R})$ with $\psi(P)=0$) is a torsion point of $E({\mathcal R})$. This, of course, fails if ${\mathcal R}$ is not assumed algebraic over ${\mathbb C}(t)$. By the way, in this general case,  for $P=(X,Y)$ with $\psi(P)=0$, the (multivalued) integral $\int_{\infty}^P\frac{dx}{y}$, as a function of $z$, is (for each of its branches) a ${\mathbb C}$-linear combination of periods-functions, so
  \begin{equation}
 \label{P}
 P=\pi\left(a_1 \cdot \int_{\gamma_1}\frac{dx}{y}+a_2\cdot \int_{\gamma_2}\frac{dx}{y}\right),\end{equation}
 where $\pi$ is the uniformization map for our family of elliptic curves (the inverse of the multivalued Abel-Jacobi map $Q\mapsto \int_{\infty}^Q\frac{dx}{y}$), $a_1,a_2\in {\mathbb C}$ and $\gamma_1,\gamma_2$ a basis for the integral homology. Of course $\pi$ is defined by elliptic functions.
 Non-torsion $P$s as above always exist and are transcendental over ${\mathbb C}(z)$ so one can refer to such a $P$ as a (degenerate) Painlev\'{e}  transcendent.  One   can  also call such a $P$  a {\it period-point} of our elliptic curve (family) since it is the composition of the uniformization map with a complex 
combination of periods.
 
   There are similar results for \ref{suffering} as follows. First it was proved in \cite{char} that, in  case $E$ has ordinary reduction mod $p$,  the set of solutions to $\psi_p(P)=0$ (which is, of course, a subgroup of $E(R)$) contains the group
 \begin{equation}
 \label{oracle}
 \bigcap_{n=1}^{\infty} p^n E(R)\end{equation}
 as a subgroup of finite index (whose index can be computed explicitly, cf. loc. cit.)
 On the other hand upon choosing a basis $(\overline{P}_n)$ of the (physical) Tate module of $E(k)$ (where $P_n\in E(R)$, the image $\overline{P}_n\in E(k)$  of $P_n$  generates $E(k)[p^n]$,  and $p\overline{P}_n=\overline{P}_{n-1}$)  the group \ref{oracle} contains the  element 
 \begin{equation}
 \label{q(E)}q(E)=\lim p^nP_n
 \end{equation}
  which can be canonically identified (via the identification of the formal group of $E$ with the formal group of ${\mathbb G}_m$) with the Serre-Tate parameter of $E$.   Finally it is easy to check  \cite{frob}
 that, if $E$ is defined over $\bZ_p$, then $q(E)$ is transcendental over $\bZ_p$ and the only solutions to $\psi_p(P)=0$ with $P\in E(\bZ_p^{ur})$ are torsion points. One of the morals to the story is that the Serre-Tate parameter $q(E)$  can be viewed, in some sense, as a {\it period-point} and plays the role of a (degenerate) Painlev\'{e} number; here by {\it degenerate} we mean ``corresponding to the degenerate case" when the $\lambda$s are $0$. For the ``non-degenerate case" the Painlev\'{e} numbers should then be viewed as deformations of the Serre-Tate parameter, and hence of the ``period-point" of $E$. 
 
 \subsection{Modular functions}
 One final example we shall discuss here arises in relation to modular forms.
 First we recall the classical (differential equation) picture; cf. \cite{mahler, ajm3, scanlon}.
 We begin by recalling  from \cite{ajm3} that there exists a rational function $\rho$ of one variable with coefficients in $\mathbb Q$ (which is entirely explicit) such that, for any $\beta\in {\mathcal R}\backslash {\mathbb C}$, the equation
 \begin{equation}
 \label{schw}
 \frac{2 (\d_z u)(\d_z^3 u)  - 3 (\d_z^2 u)^2}{4(\d_z u)^2}+(\d_z u)^2\cdot \rho(u)=\beta
 \end{equation}
  is ``constant on isogeny classes" in the following sense: if $u\in {\mathcal R}\backslash {\mathbb C}$,  is a solution to  \ref{schw} then any $u^*\in {\mathcal R}\backslash {\mathbb C}$ in the same isogeny class as $u$ is also a solution to \ref{schw}. Here two elements  $u_1,u_2 \in {\mathcal R}$ are said to be in the same isogeny class if the two elliptic curves over ${\mathcal R}$ with $j$-invariants $u_1$ and $u_2$ respectively are isogenous. (Of course, the first term in \ref{schw} is the classical Schwartzian.) 
 For work on the algebraic relations among the solutions of \ref{schw}, leading to analogues for the $j$ function of the Ax-Lindemann-Weiestrass theorem, we refer to \cite{pilaa, pila, scanlon}; the  solutions of \ref{schw} one looks at  are of the form
 \begin{equation}
 \label{j}
 j\left(\frac{a \tau(z)+b}{c\tau(z)+d}\right)
 \end{equation}
 with $\tau(z)$ an algebraic function, and $a,b,c,d\in {\mathbb C}$; \ref{j} should be viewed as an analogue of \ref{exponential} and \ref{P} and hence as a ``period-point", as it is the composition of the uniformization map $j$ with a complex Mobius transformation applied to a ``periods function" $\tau(z)$ (of some elliptic curve).
 So the function $j$ in \ref{j} is analogous to the function $\pi$ in \ref{P} and to the function $\exp$ in \ref{exponential}.
 The constants $a,b,c,d$ in \ref{j} are analogues of the constants $a_1,a_2$ in \ref{P} and are also analogues of the constants appearing as coefficients of $\beta$ in \ref{exponential}, if $\beta$ in \ref{exponential} is, say, a polynomial in $z$. Finally the function $\tau(z)$ in \ref{j} is analogous to the integrals in \ref{P} and to the integral in \ref{exponential}.
 
 The above story has an arithmetic analogue as follows; cf. \cite{difmod, Barcau, book}. Let $X_1(N)$ be the  complete modular curve over $R$ of level $N>4$
and let $L_{X_1(N)}$ be the line bundle on $X_1(N)$ with the property that
the sections of its various powers are the classical modular forms
on $\Gamma_1(N)$ of various weights. Let $Y_1(N)$ be $X_1(N)$ from which one removes the 
cusps. Let
$L_{Y_1(N)}\ra Y_1(N)$ be the restriction of the above line bundle and let $V$ be
$L_{Y_1(N)}$ with the zero section removed. The $R$-points of $V$ correspond then to certain Weierstrass cubics equipped with level structure. Let $M^r$ be the ring of $\d$-functions of order $r$ on $V$. For $w\in \bZ[\phi_p]\subset End_{\bZ-\text{mod}}(R)$ we denote by $M^r(w)$ the space of all members of $M^r$ that have ``weight" $w$ (i.e. on which ${\mathbb G}_m$ acts via the character $\lambda\mapsto \lambda^w$); the elements of $M^r(w)$ are called {\it $\d$-modular forms} of order $r$ and weight $w$. The theory in \cite{book} provides some remarkable $\d$-modular forms $f_p^r\in M^r(-1-\phi^r)$, $r\geq 1$, 
which have a property (called {\it isogeny covariance}) that is stronger than the property of being Hecke eigenforms; cf. \cite{difmod, book}. This property    has no analogue in the classical theory and we refer to loc.cit. for its definition. Then, for any $\beta\in R$, one can consider the arithmetic differential equation of order $3$ of the following form
\begin{equation}
\label{new}
\frac{f_p^2(P)^{1+\phi_p}}{f_p^1(P)^{1+\phi_p^2}}=\beta,
\end{equation}
with  solutions $P\in V(R)$. These equations  are ``constant on isogeny classes": if $P$ is a solution to \ref{new} (with the denominator in \ref{new} invertible)
then any $P^*\in V(R)$ in the same isogeny class as $P$ is also a solution to \ref{new}; here two points of $V$ are said to be in the same isogeny class if and only if there is an isogeny between the Weierstrass cubics representing the two points which is compatible with the standard $1$-forms on the cubics (but not necessarily with the level structures). 
Clearly \ref{new} can be viewed as an analogue of \ref{schw}.  Note that both \ref{schw} and \ref{new} have order $3$ which seems to be expected, since the dimension of $SL_2$ is $3$ and $SL_2$ is the group behind the isogeny equivalence relation. 
It is therefore all the more surprising to see, as we shall explain in what follows, that there is an order $2$ equation that is analogous to \ref{schw}. Indeed, let $V_{ord}$, $M_{ord}(w)$ be the spaces obtained by replacing, everywhere in the construction above, the curve $Y_1(N)$ with its ordinary locus (i.e. the locus where the Eisenstein series $E_{p-1}$ is invertible). Then an entirely new actor joins the scene: a remarkable isogeny covariant $\d$-modular form of order $1$ and weight $\phi-1$, introduced in \cite{Barcau}, and denoted by $f_p^{\partial}\in M^1_{ord}(\phi-1)$. One can then consider, for any $\beta\in R$,  the arithmetic differential equation of order $2$
\begin{equation}
\label{brand}
f_p^1(P)^{\phi_p-1}\cdot f_p^{\partial}(P)^{\phi_p+1}=\alpha,
\end{equation}
with  solutions $P\in V_{ord}(R)$. This equation is, again, ``constant on isogeny classes": if $P$ is one of its  solutions (with $E_{p-1}(P)$ invertible)
then any $P^*\in V_{ord}(R)$ in the same isogeny class as $P$ is also one of its solutions. 
In fact \ref{new} is a consequence of \ref{brand} for 
$\beta=\alpha+p(1+\alpha^{1-\phi})+p^2\alpha^{-\phi^2}$.
The fact that \ref{brand} has order $2$ rather than $3$ may seem  somewhat mysterious but here is one way to think about this phenomenon.
It is a fact (known to Mahler \cite{mahler}, see also \cite{ajm3}, p.37)  that there is no order $2$ algebraic differential equation that is ``constant on isogeny classes". This is equivalent to the fact that if $j$ is the classical holomorphic $j$-function, $j=j(z)$, defined on the upper half plane then $j, \d_z j, \d_z^2 j$ are algebraically independent over ${\mathbb C}$. 
Now these series belong to $\bZ((q))=\bZ[[q]][q^{-1}]$, $q=\exp(2 \pi i z)$. So setting 
$\theta=q\frac{d}{dq}$ the series $j, \theta j, \theta^2 j$ are algebraically independent
over ${\mathbb Q}$ as elements of $\bZ((q))$. However it turns out that
the images of the series $j, \theta j$ in ${\mathbb F}_p((q))$ are algebraically dependent over ${\mathbb F}_p$; cf. \cite{lang}; so 
$j(q)$ mod $p$ satisfies an algebraic differential equation of order $1$!
This ``drop in order" from $3$ to $1$ is a remarkable phenomenon in characteristic $p$. On the other hand we have a ``drop in order" from $3$ to $2$ between \ref{schw} and \ref{brand} in the framework of $\d$-modular forms and it is possible to prove that, in some precise sense, no ``drop in order" from $3$ to $1$ exists in this framework. So the drop from $3$ to $2$ can be viewed as a partial lift to characteristic zero, in  $\d$-geometry, of the drop from $3$ to $1$ that exists in characteristic $p$. 
  The existence of such a lift is surprising but
 it is just one instance of a more general phenomenon where characteristic $p$ algebro-geometric objects that do not lift to characteristic zero in algebraic geometry {\it do} lift to characteristic zero in $\d$-geometry; cf. \cite{book} for a series of examples of this.
 
 By analogy with Pila's results \cite{pila} it would be interesting to prove that there are no algebraic relations among the solutions of \ref{new} (or \ref{brand}) except those given by ``isogeny". The simplest case of this would be the following conjecture (which pertains to \ref{brand}). {\it Let $P,P^*\in V_{ord}(R)$
 be two points corresponding to Weierstrass elliptic curves with coefficients in a number field in which $p$ splits completely. Assume $f^{\partial}(P)=f^{\partial}(P^*)$. Then  the two elliptic curves are isogenous.}

Note that there is a cocompact analogue of both \ref{schw} and
 \ref{new}, \ref{brand} (cf. \cite{ullmo} and \cite{book} respectively) with transcendence results in the functional case \cite{ullmo}.

\section{Arithmetic differential Galois groups and periods}

Galois groups are a measure of the  relations among solutions of  equations.  Galois' original work is about algebraic equations and algebraic relations among their roots.  Picard-Vessiot theory (or more generally Kolchin theory \cite{kolchin}) is about linear differential equations (respectively differential equations arising from Maurer-Cartan connections on algebraic groups) and about algebraic relations among their solutions. 
More generally  \cite{land} considers differential equations with parameters  and  differential algebraic relations among solutions; cf. also \cite{CS} for the  case of linear equations.
The same idea
can be applied, in certain cases,  to arithmetic differential equations and the $\d$-algebraic relations among their solutions; this can be done, e.g.,  for equations of the form \ref{lineara} (cf. \cite{adel3}), or 
\ref{difference} (cf. \cite{book}), or
\ref{suffering}, at least in the degenerate case when the $\lambda$s are $0$ (cf. \cite{book}).  No general Galois theory of arithmetic differential equations is yet available. However, if such a theory can be developed and, at the same time, various ``$p$-adic periods" appearing in algebraic geometry  can be shown to satisfy ``nice" arithmetic differential equations, then it is conceivable that one could systematically understand relations among ``$p$-adic periods" (and possibly motivic Galois groups \cite{andre})  through a Galois theory of arithmetic differential equations. The  possibility suggested above is compatible with suggestions made by Manin in \cite{manin}. 
Whether or not the ``$p$-adic periods" referred to above should be related to the ones appearing in the standard comparison isomorphisms \cite{andre}  is not clear at this time. Let us close by making some comments on this issue, partly inspired by comments  in \cite{manin}. In the functional case \ref{expo}, \ref{P}, \ref{j},
the solutions of our differential equations, which we would like to think of as {\it period points} on various schemes,  are compositions of ``uniformizing functions"
$\exp,\pi,j$ with complex linear (respectively ``Mobius") combinations of functions that have the flavor of {\it  periods functions}. So the period points are not really periods functions but rather images of such via uniformization maps. Such uniformization maps {\it do not exist a priori} in the $p$-adic situation (e.g. for elliptic curves no uniformization is known for the good reduction case!) Therefore the solutions to  \ref{e}, \ref{suffering}, \ref{new} should be viewed as analogues of {\it period points} and only indirectly as analogues of {\it $p$-adic periods} in the motivic setting \cite{andre}. As for the ``motivic periods"  it is worth recalling from \cite{andre} the general philosophy as follows. Say
that $X$ is a smooth projective variety over a field $k$ of characteristic $0$. For $k\subset {\mathbb C}$ one has a deRham versus Betti comparison isomorphism $H_{dR}(X)\otimes_k{\mathbb C}\simeq 
H_{Betti}(X)\otimes_{\mathbb Q}{\mathbb C}$, due to Grothendieck, where $H_{dR}$ is  the algebraic deRham cohomology of $X/k$ and $H_{Betti}$ is the Betti cohomology with ${\mathbb Q}$-coefficients; this isomorphism, expressed in bases of $H_{dR}(X)$ and $H_{Betti}(X)$ provides a complex matrix whose entries are called the {\it complex periods} of $X$. 
By the way one also has in this case the $p$-adic \'{e}tale versus Betti comparison isomorphism  $H_{et}(X)\simeq 
H_{Betti}(X)\otimes_{\mathbb Q}{\mathbb Q}_p$, due to Grothendieck and Artin.
Assume until further notice that $k$ is a number field. For any $k$-subvariety $Z$ of any power $X^n$ the comparison isomorphism for $X^n$ sends the class of $Z$ into the corresponding class of $Z$; via K\"{u}nneth, this induces an algebraic relation over $k$ between the complex periods. Grothendieck's period conjecture says roughly that {\it all} the relations among periods ``come from" $Z$s as above. Also one defines the {\it motivic Galois group} of $X$ as, roughly, the algebraic ${\mathbb Q}$-subgroup of $GL(H_{Betti}(X))$ that fixes all the classes of $Z$s as above (viewed, via K\"{u}nneth,  as tensors on $H_{Betti}(X)$). 
In the $p$-adic case one has a similar picture. Indeed for $k$ a finite extension of ${\mathbb Q}_p$ one has the deRham versus $p$-adic \'{e}tale comparison isomorphism  $H_{dR}(X)\otimes_k B_{dR}\simeq 
H_{et}(X)\otimes_{{\mathbb Q}_p}B_{dR}$, due to Faltings, where $B_{dR}$ is Fontaine's field of $p$-adic periods \cite{fontaine}.
The latter yields, again, period matrices with entries in $B_{dR}$; these {\it $p$-adic periods}, lying in $B_{dR}$,  are a priori rather different from the {\it period points}, referred to above, which are $R$-points (hence with coordinates lying in $R$) of certain $R$-schemes (the hyperbola $H={\mathbb G}_m$ over $R$, elliptic curves $E$ over $R$, ${\mathbb G}_m$-bundles $V$ over modular curves over $R$, etc.) A link between these two types of periods would allow one to consider arithmetic differential equations in the context of periods lying in $B_{dR}$. We end by recalling the definition of $B_{dR}$ \cite{fontaine} and the $p$-adic periods in \cite{colmez}. 
Assume, for simplicity, that $k$ is unramified over ${\mathbb Q}_p$.
One lets $\cO$ be the ring of integers of ${\mathbb C}_p:=\widehat{\overline{\mathbb Q}}_p$, one lets $R_{\cO}$ be the projective limit of $\cO/p\cO\leftarrow \cO/p\cO\leftarrow...$ where the arrows are the $p$-power Frobenius, 
one considers the canonical map $\theta:W(R_{\cO})\ra \cO$ (where $W$ stands for the $p$-typical Witt functor), one lets $A_{inf}$ be the completion of $W(R_{\cO})$ with respect to $(p,Ker(\theta))$, one lets 
$\theta_k:A_{inf}[1/p]\ra {\mathbb C}_p$ be  the induced map, one lets $B^+_{dR}$ be the completion of $A_{inf}[1/p]$ with respect to $Ker(\theta_k)$, which is a complete discrete valuation ring with residue field ${\mathbb C}_p$, and  one  lets $B_{dR}=Frac(B^+_{dR})$. Uniformizers in $B_{dR}$ can be constructed from any basis of the Tate module of ${\mathbb G}_m({\mathbb C}_p)$ as follows. Let $\epsilon_n\in {\mathbb C}_p$ be primitive $p^n$-roots of $1$, $\epsilon_n^p=\epsilon_{n-1}$, and let $\tilde{\epsilon}_n\in A_{inf}$ be lifts of $\epsilon_n$; then $\tilde{\epsilon}=\lim \tilde{\epsilon}_n^{p^n}$
exists in $A_{inf}$ and $\log_p \tilde{\epsilon}\in B^+_{dR}$ is a uniformizer of $B^+_{dR}$;
note that the limit $\tilde{\epsilon}$ is reminiscent of  \ref{q(E)}. Furthermore, for any abelian variety $A$ over $k$, any $1$-form $\omega$ on $A$ and any element $\gamma$ in the Tate module of $A({\mathbb C}_p)$ there is, by work of Colmez \cite{colmez},  a well defined period $\int_{\gamma}\omega=\lim p^n I_n \in B^+_{dR}$, where $I_n$ is a certain integral on $A(B_{dR})$. This  construction is, again, reminiscent of \ref{q(E)}; note however that the image of $\int_{\gamma}\omega$ in ${\mathbb C}_p$ turns out to be $0$ so (for $E$ coming from an integral model of $A/k$) this  image does not coincide with  $q(E)\in R\subset {\mathbb C}_p$; indeed    $q(E)\neq 0$ for  $E$ without complex multiplication.

\end{document}